\documentclass[12pt]{article}

\usepackage{amssymb,amsmath,amsthm,
mathabx,MnSymbol}

\usepackage[all]{xy}
\usepackage[active]{srcltx}

\theoremstyle{definition}
\newtheorem*{definition}{DEFINITIONS}
\newtheorem*{fac}{FACT}
\newtheorem*{fact}{FACTS}
\newtheorem*{example}{EXAMPLES}
\newtheorem*{exampl}{EXAMPLE}

\renewcommand{\le}{\leqslant}

\newcommand{\ci}{
\begin{picture}(6,6)
\put(3,3){\circle*{3}}
\end{picture}}

\newcommand{\mrof}%
{\unitlength 1.2pt
\begin{picture}(8,0)(7,13.1)
\put(10,15){\oval(4,4)[l]}
\qbezier(10,17)(12,17)(14,16)
\qbezier(10,13)(12,13)(14,14)
\end{picture}}

\newcommand{\form}%
{\unitlength 1.2pt
\begin{picture}(8,0)(5,12.6)
\put(10,15){\oval(4,4)[r]}
\qbezier(10,17)(8,17)(6,16)
\qbezier(10,13)(8,13)(6,14)
\end{picture}}

\newcommand{\fform}%
{\unitlength 2pt
\begin{picture}(8,0)(5,14)
\put(10,15){\oval(4,4)[r]}
\qbezier(10,17)(8,17)(6,16)
\qbezier(10,13)(8,13)(6,14)
\end{picture}}

\newcommand{\lin}{\,\frac{}{\ \quad }\,}

\DeclareMathOperator{\Ker}{Ker}
\DeclareMathOperator{\End}{End}
\DeclareMathOperator{\ind}{ind}
\DeclareMathOperator{\diag}{diag}

\begin{document}

\title{Representations of quivers and
mixed graphs\thanks{This is a preprint
of Chapter 34 from \emph{Handbook of
Linear Algebra} (Edited by L. Hogben),
Second Edition, Chapman \& Hall/CRC,
Boca Raton, FL, 2014.}}

\author{Roger A. Horn\\The
University of Utah, USA \and Vladimir
V. Sergeichuk\\ Institute of
Mathematics, Kiev, Ukraine}

\date{}
 \maketitle

\begin{abstract}
This is a survey article for
\emph{Handbook of Linear Algebra}, 2nd
ed., Chapman \& Hall/CRC, 2014. An
informal introduction to
representations of quivers and finite
dimensional algebras from a linear
algebraist's point of view is given.
The notion of quiver representations is
extended to representations of mixed
graphs, which permits one to study
systems of linear mappings and bilinear
or sesquilinear forms. The problem of
classifying such systems is reduced to
the problem of classifying systems of
linear mappings.

\end{abstract}

\section*{Introduction}
\noindent In Sections
\ref{ser1}--\ref{ser3}, we give an
informal introduction to quivers from a
linear algebraist's point of view.
Exact definitions, results, and their
proofs can be found in surveys
\cite{weyman,physenc} and monographs
\cite{simson,ausl,dl+rin,gab_roi,haz-kir,%
pierce,rin}.

After Gabriel's article \cite{gab}, in
which the notions of a quiver and its
representations were introduced, it
became clear that a whole range of
problems about systems of linear
mappings can be formulated and studied
in a uniform way. Quivers arise
naturally in many areas of mathematics
(representation theory, algebraic and
differential geometry, number theory,
Kac-Moody algebras, quantum groups,
geometric invariant theory) and physics
(string theory, supersymmetry, black
holes, particle physics). Each finite
dimensional algebra can be given by a
quiver with relations, and
representations of the algebra can be
identified with representations of this
quiver; that is, with finite systems of
linear mappings satisfying some
relations. Thus, the modern theory of
representations of finite dimensional
algebras can be considered as a branch
of linear algebra.

In Sections \ref{ser4} and \ref{ser5},
we extend the notion of quiver
representations to representations of
mixed graphs, which permits one to
study systems of linear mappings and
bilinear or sesquilinear forms. We
reduce the problem of classifying such
systems to the problem of classifying
systems of linear mappings.


\section{Systems of linear mappings as
representations of quivers}\label{ser1}

\begin{definition}
\noindent A {\bf quiver} $Q$ is a
directed graph where multiple loops and
multiple arrows between two vertices
are allowed. We suppose that the
vertices of $Q$ are $1,\dots,t$ and
denote by $\alpha :i\longrightarrow j$
an arrow $\alpha $ from a vertex $i$ to
a vertex $j$.

A {\bf representation} $\cal A$ of $Q$
over a field $\mathbb F$ is given by
assigning to each vertex $i$ a finite
dimensional vector space $A_i$ over
$\mathbb F$ and to each arrow
$\alpha:i\longrightarrow j$ a linear
mapping $A_{\alpha}:A_i\to A_j$.

The {\bf dimension} of  $\cal A$ is the
vector $${\bf z}=(\dim A_1,\dots, \dim
A_t).$$

A {\bf morphism} $\varphi: {\cal A}\to
{\cal B}$ between representations
${\cal A}$ and ${\cal B}$ of $Q$ is a
family of linear mappings
$$\varphi_1:A_1\to B_1,\ \dots,\
\varphi_t:A_t\to B_t$$ such that the
diagram
\[
\xymatrix{
A_i\ar[r]^{A_{\alpha} }
\ar[d]_{\varphi_i}
&A_j\ar[d]^{\varphi_j}
\\
B_i\ar[r]^{B_{\alpha} }&B_j}
\]
is commutative (i.e.,
$\varphi_jA_{\alpha }=B_{\alpha
}\varphi _i$) for each arrow $\alpha
:i\longrightarrow j$.

An {\bf isomorphism} $\varphi: {\cal A}
\
\begin{matrix}\sim\; \\[-9pt] \to\\[-9pt]{}
\end{matrix}\; {\cal B}$
is a morphism $\varphi: {\cal A}\to
{\cal B}$ in which all $\varphi _i$ are
bijections.

The {\bf direct sum} $ {\cal A}\oplus
{\cal B}$
 of representations
${\cal A}$ and ${\cal B}$ of $Q$ is the
representation of $Q$ defined by \[
(A\oplus B)_i:=A_i\oplus B_i,\qquad
(A\oplus B)_{\alpha }:=A_{\alpha
}\oplus B_{\alpha }\] for all vertices
$i$ and arrows $\alpha $. (The direct
sum of linear mappings $A:U\to V$ and
$A':U'\to V'$ is the linear mapping
$$A\oplus A':U\oplus U'\to V\oplus V'$$
defined by $$A\oplus A':u+u'\mapsto
Au+A'u'\qquad\text{for all $u\in U$ and
$u'\in U'$}.)$$

A representation of nonzero dimension
is {\bf indecomposable} if it is not
isomorphic to a direct sum of
representations of smaller dimensions.
\end{definition}

\begin{fac}
\begin{enumerate}
\enlargethispage{12pt}
 \item[] 
\emph{The Krull--Schmidt theorem}
\cite[Corollary 2.4.2]{haz-kir}:
Each representation of a quiver is
isomorphic to a direct sum of
indecomposable representations.
This direct sum is uniquely
determined, up to permutation and
isomorphisms of direct summands;
that is, if \[ {\cal
A}_1\oplus\dots\oplus {\cal
A}_r\simeq {\mathcal
B}_1\oplus\dots\oplus \mathcal
B_s,\] in which all $\mathcal A_i$
and $\mathcal B_j$ are
indecomposable representations,
then $r=s$ and all $\mathcal
A_i\simeq \mathcal B_i$ after a
suitable renumbering of $\mathcal
A_1,\dots,\mathcal A_r$.
\end{enumerate}
 \end{fac}

\begin{example}
\begin{enumerate}
\item \label{jjgt} Each
    representation
    \[\raisebox{20pt}
     
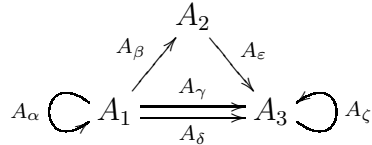
{\xymatrix@R=20pt@C=10pt{
 &{A_2}&\\
 {A_1} \save !<-2mm,0cm>
\ar@(ul,dl)@{->}_{A_{\alpha }}
\restore
 \ar@{->}[ur]^{A_{\beta }}
 \ar@<0.4ex>[rr]^{A_{\gamma }}
 \ar@<-0.4ex>[rr]_{A_{\delta }}
 &&{A_3}
 \ar@{<-}[ul]_{A_{\varepsilon }}
 \save !<2mm,0cm>
\ar@{<-}@(ur,dr)^{A_{\zeta}}
\restore }}
\]
of the quiver
 \[\raisebox{20pt}
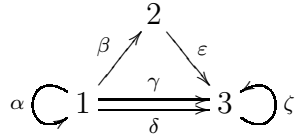
{\xymatrix@R=20pt@C=15pt{
 &{2}&\\
 {1} \ar@(ul,dl)@{->}_{\alpha }
 \ar@{->}[ur]^{\beta }
 \ar@<0.4ex>[rr]^{\gamma }
 \ar@<-0.4ex>@{->}[rr]_{\delta }
 &&{3} \ar@{<-}[ul]_{\varepsilon }
\ar@{<-}@(ur,dr)^{\zeta} }} \] over
a field $\mathbb F$ is a system of
vector spaces $A_1,A_2,A_3$ over
$\mathbb F$ and linear mappings
$A_{\alpha }:A_1\to A_1$, $A_{\beta
}:A_1\to A_2$,\dots

\item Consider the problems of
    classifying representations of
    the quivers
\begin{equation*}\label{sol}
\xymatrix{
 {1} \ar[r]&{2}}
,\quad
\xymatrix{
 {1} \ar@(dr,ur)@{->}}\quad\ \ \ ,\quad
  \xymatrix{
 {1}
 \ar@<0.4ex>[r]
 \ar@<-0.4ex>[r]
 &{2}},\quad
  \xymatrix{
 {1}\ar@<0.4ex>[r]&{2}
 \ar@<0.4ex>[l]
 },
 \quad\quad\quad
\xymatrix{
 {1} \ar@(ul,dl)@{->} \ar@(dr,ur)@{->}}\quad\
\end{equation*}

\begin{itemize}
 \item Each matrix $A\in
     \mathbb F^{m\times n}$
     defines the representation
     $\mathbb F^n\xrightarrow{\
     A \ }\mathbb F^m$ of the
     quiver $
     {1}\longrightarrow{2}$ by
     assigning to its arrow the
     linear mapping $x\mapsto
     Ax$ with $x\in \mathbb
     F^n$. Thus,
  the problem of classifying
     representations of the
     quiver $
     {1}\longrightarrow{2}$ is
     the canonical form problem
     for matrices under
     equivalence
     transformations $A\mapsto
     R^{-1}AS$ with nonsingular
     $R$ and $S$. Its canonical
     matrices are $I\oplus 0$,
     and so each representation
     is isomorphic to a direct
     sum, uniquely determined
     up to permutations of
     summands, of
     representations of the
     form
\begin{equation}\label{key}
\mathbb F \xrightarrow{\ I_1\ }\mathbb F,\qquad
0 \xrightarrow{\ 0_{10}\ }\mathbb F,\qquad \mathbb F \xrightarrow{\ 0_{01}\ }0
\end{equation}
(it is agreed that $\mathbb
F^0=0$ and there exist exactly
one matrix $0_{n0}$ of size
$n\times 0$ and exactly one
matrix $0_{0n}$ of size
$0\times n$ for every
nonnegative integer $n$; they
are the matrices of linear
mappings $0\to \mathbb F^n$ and
$\mathbb F^n\to 0$).

 \item\label{ex2} The problem
     of classifying
     representations of the
     quiver
     $1\!\!\righttoleftarrow$
     is the canonical form
     problem for an $m\times m$
     matrix $A$ over a field
     $\mathbb F$ under
     similarity transformations
     $S^{-1}AS$ with
     nonsingular $S\in \mathbb
     F^{m\times m}$. Its
     canonical matrix is a
     direct sum of companion
     matrices
\begin{equation}\label{kut}
C_n(q)
      =
\left[\begin{array}{r@{\quad}r@{\quad}r@{\quad}c} 0&&
0&-c_n\\1&\ddots&&\vdots
\\&\ddots&0&-c_2\\
0&&1&-c_1 \end{array}\right]
\end{equation}
whose characteristic
polynomials
$$q(x)=x^n+c_1x^{n-1}+\dots+c_n$$
are powers of irreducible
polynomials. This canonical
matrix is called the elementary
divisors rational canonical
form of $A$, or the Frobenius
canonical form of $A$. Thus,
each representation of
$1\!\!\righttoleftarrow$ is
isomorphic to a direct sum,
uniquely determined up to
permutation of summands, of
representations of the form ${
 {\mathbb
 F^n}\!\righttoleftarrow
 {\!\scriptstyle C_n(q)}.}$ If
$\mathbb F$ is an algebraically
closed field, then a Jordan
block $J_n(\lambda )$ can be
taken instead of $C_n(q)$.
\medskip

 \item\label{ex3} The problem
     of classifying
     representations of the
     quiver $
     \!\!\xymatrix@=15pt{{1}
     \ar@<0.4ex>[r]
 \ar@<-0.4ex>[r] &{2}}\!$ is
 the canonical form problem for
 pairs $(A,B)$ of matrices of
 the same size under
 equivalence transformations $$
 (R^{-1}AS,R^{-1}AS)\quad\text{with
 nonsingular $R$ and $S$.}$$ By
Kronecker's theorem on pencils
of matrices, each
representation of $
     \!\!\xymatrix@=15pt{{1}
     \ar@<0.4ex>[r]
 \ar@<-0.4ex>[r] &{2}}\!$ is
isomorphic to a direct sum,
uniquely determined up to
permutations of summands, of
representations of the form
\begin{equation}\label{lic}
 \xymatrix{
 {\mathbb F^n}
 \ar@<0.4ex>[r]^{I_n}
 \ar@<-0.4ex>[r]_{C_n(q)}
 &{\mathbb F^n}},
              \
              \xymatrix{
 {\mathbb F^n}
 \ar@<0.4ex>[r]^{J_n(0)}
 \ar@<-0.4ex>[r]_{I_n}
 &{\mathbb F^n}},
                \
 \xymatrix{
 {\mathbb F^{n}}
 \ar@<0.4ex>[r]^{L_n}
 \ar@<-0.4ex>[r]_{R_n}
 &{\mathbb F^{n-1}}},
                \
 \xymatrix{
 {\mathbb F^{n-1}}
 \ar@<0.4ex>[r]^{\ \ L_n^T}
 \ar@<-0.4ex>[r]_{\ \ R_n^T}
 &{\mathbb F^{n}}},
\end{equation}
in which $n=1,2,\dots$,
\begin{equation}\label{1.4}
L_n:=\begin{bmatrix}
1&0&&0\\&\ddots&\ddots&\\0&&1&0
\end{bmatrix},\quad
R_n:=\begin{bmatrix}
0&1&&0\\&\ddots&\ddots&\\0&&0&1
\end{bmatrix}\quad\text{($(n-1)$-by-$n$)},
\end{equation} $L_1=R_1=0_{01}$, and $C_n(q)$ is a
block \eqref{kut}, which can be
replaced by a Jordan block if
$\mathbb F$ is algebraically
closed.

 \item\label{ex4} The problem of
     classifying representations of
     the quiver $
  \!\!\xymatrix@=15pt{
 {1}\ar@<0.4ex>[r]&{2}
 \ar@<0.4ex>[l] }\!$ is the
 canonical form problem for
 pairs $(A,B)$ of $p\times q$
 and $q\times p$ matrices under
 contragredient equivalence
 transformations
 $$(R^{-1}AS,S^{-1}AR)\quad\text{with
 nonsingular $R$ and $S$.}$$
Dobrovol$'$skaya and Ponomarev
\cite{dob+pon} (see also
\cite{hor+mer}) proved that
each representation of $
  \!\!\xymatrix@=15pt{
 {1}\ar@<0.4ex>[r]&{2}
 \ar@<0.4ex>[l] }\!$ is
isomorphic to a direct sum,
determined uniquely up to
permutation of summands, of
representations of the form
\[
 \xymatrix{
 {\mathbb F^n}
 \ar@<0.4ex>[r]^{I_n}
 &{\mathbb F^n}
 \ar@<0.4ex>[l]^{C_n(q)}
 },
              \quad
 \xymatrix{
 {\mathbb F^n}
 \ar@<0.4ex>[r]^{J_n(0)}
 &{\mathbb F^n}\ar@<0.4ex>[l]^{I_n}
 },
               \quad
 \xymatrix{
 {\mathbb F^{n}}
 \ar@<0.4ex>[r]^{L_n}
 &{\mathbb F^{n-1}}
 \ar@<0.4ex>[l]^{R_n^T}},
               \quad
 \xymatrix{
 {\mathbb F^{n-1}}
 \ar@<0.4ex>[r]^{\ \ L_n^T}
 &{\mathbb F^{n}}
 \ar@<0.4ex>[l]^{\ \ \ R_n}},
\]
in which $n=1,2,\dots$, the
matrices $L_n$ and $R_n$ are
defined in \eqref{1.4}, and
$C_n(q)$ is a block
\eqref{kut}, which can be
replaced by a Jordan block if
$\mathbb F$ is algebraically
closed.
\end{itemize}
 \end{enumerate}
\end{example}


\section{Tame and wild quivers}\label{ser2}

The problem of classifying pairs of
$n\times n$ matrices up to similarity
transformations
\begin{equation*}\label{gpw}
(A,B)\mapsto (S^{-1}AS,S^{-1}BS)\qquad
\text{with nonsingular } S
\end{equation*}
(i.e., representations of the quiver
 $\lefttorightarrow  \!\!1\!\!\righttoleftarrow$)
plays a special role in the theory of
quiver representations: it contains the
problem of classifying representations
of each quiver.

\begin{definition}

\noindent A quiver is of {\bf finite
type} if it has only finitely many
nonisomorphic indecomposable
representations. A quiver is of {\bf
wild type} if the problem of
classifying its representations
contains the problem of classifying
matrix pairs up to similarity,
otherwise the quiver is of {\bf tame
type} (see formal definitions in
\cite[Section 14.10]{gab_roi}).

The {\bf Tits quadratic form}
$q_Q:\mathbb Z^t\to \mathbb Z$ of a
quiver $Q$ with vertices $1,\dots,t$ is
the form
\begin{equation}\label{kur}
q_Q(x_1,\dots,x_t):=x_1^2+\dots+x_t^2-
\sum_{i\longrightarrow j}x_ix_j
\end{equation}
in which the sum is taken over all
arrows of the quiver.

\end{definition}

\begin{fact}
\begin{enumerate}

 \item 
The problem of classifying pairs of
commuting nilpotent matrices up to
similarity contains the problem of
classifying arbitrary matrix pairs
up to similarity (see
\cite{gel-pon} and Example 1).

 \item 
The problem of classifying matrix
pairs up to similarity contains the
problem of classifying
representations of any quiver (see
\cite{gel-pon,bel-ser_compl} and
Example 2).

 \item 
 \emph{Gabriel's theorem}
\cite{gab}: Let $Q$ be a connected
quiver with $t$ vertices.
\begin{itemize}
  \item $Q$ is of finite type
      if and only if the Tits
      form $q_Q$ (considered as
      a form over $\mathbb R$)
      is positive definite,  if
      and only if $Q$ can be
      obtained by directing
      edges in one of the
      Dynkin diagrams
\begin{equation}\label{dfe}
\begin{aligned}
&A_t\ \xymatrix@=10pt@R=0,5pt{
&&&&\\
*{\ci}\ar@{-}[r]&*{\ci}\ar@{-}[r]&
*{{\ci}\
\cdots\ {\ci}}
\ar@{-}[r]&*{\ci}\ar@{-}[r] &*{\ci}}
\qquad
&&
D_t\ \xymatrix@=10pt@R=0,5pt{&&&&*{\ci}\\
*{
\ \ci}\ar@{-}[r]&*{\ci}\ar@{-}[r]&
*{{\ci}\
\cdots\ {\ci}}
\ar@{-}[r]&*{\ci}\ar@{-}[ur]
 \ar@{-}[dr]&\\ &&&&*{\ci}}
          \\[-5pt]
&
\begin{matrix}
 \\
E_6
\end{matrix}
\ \xymatrix@=10pt{
&&*{\ci}\ar@{-}[d]&&\\
*{\ci}\ar@{-}[r]&*{\ci}\ar@{-}[r]&
*{\ci}
\ar@{-}[r]&*{\ci}\ar@{-}[r] &*{\ci}}
&&
\begin{matrix}
 \\
E_7
\end{matrix}
\ \xymatrix@=10pt{
&&*{\ci}\ar@{-}[d]&&\\
*{\ci}\ar@{-}[r]&*{\ci}\ar@{-}[r]&
*{\ci}
\ar@{-}[r]&*{\ci}\ar@{-}[r] &*{\ci}\ar@{-}[r] &*{\ci}}
          \\
&\begin{matrix}
 \\
E_8
\end{matrix}
\ \xymatrix@=10pt{
&&*{\ci}\ar@{-}[d]&&\\
*{\ci}\ar@{-}[r]&*{\ci}\ar@{-}[r]&
*{\ci}
\ar@{-}[r]&*{\ci}\ar@{-}[r] &*{\ci}\ar@{-}[r] &*{\ci}\ar@{-}[r] &*{\ci}}
\end{aligned}
\end{equation}

(the index is the number of
vertices).

  \item Let $Q$ be of finite
      type and let ${\bf
      z}=(z_1,\dots,z_t)$ be an
      integer vector with
      nonnegative components.
      There exists an
      indecomposable
      representation of
      dimension $\bf z$ if and
      only if $q_Q({\bf z})=1$;
      this representation is
      determined by $\bf z$
      uniquely up to
      isomorphism.
      (Representations of
      quivers of finite type
      were classified by
      Gabriel \cite{gab}; see
      also \cite[Theorem
      2.6.1]{haz-kir}.)
\end{itemize}

 \item 
 \emph{The
 Donovan--Freislich--Nazarova
 theorem} \cite{don1,naz}:  Let $Q$
 be a connected quiver with $t$
 vertices.
 \begin{itemize}
 \item $Q$ is of tame type if
     and only if the Tits form
     $q_Q$ is positive
     semidefinite, if and only
     if $Q$ can be obtained by
     directing edges in one of
     the Dynkin diagrams
     \eqref{dfe} or the
     extended Dynkin diagrams
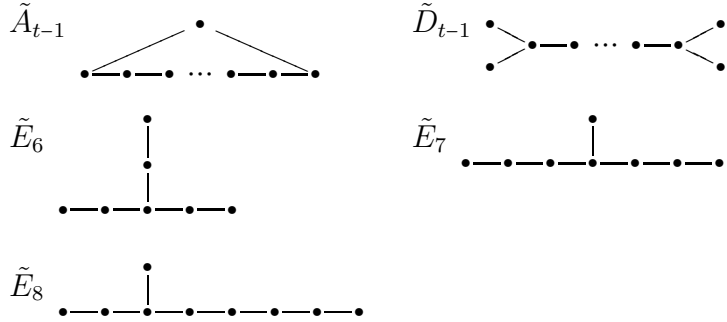
\begin{equation*}\label{dfse}
\begin{aligned}
&\tilde A_{t-1}\ \xymatrix@=10pt{
&&{\ci}\ar@{-}[lld]\ar@{-}[rrd]&&\\
*{\ci}\ar@{-}[r]&*{\ci}\ar@{-}[r]&
*{{\ci}\
\cdots\ {\ci}}
\ar@{-}[r]&*{\ci}\ar@{-}[r] &*{\ci}}
\qquad
&&
\tilde D_{t-1}\ \xymatrix@=10pt@R=0,5pt{*{\ci}&&&&*{\ci}\\
&*{\ci}\ar@{-}[r]\ar@{-}[lu]\ar@{-}[ld]&
*{{\ci}\
\cdots\ {\ci}}
\ar@{-}[r]&*{\ci}\ar@{-}[ur]
 \ar@{-}[dr]&\\ *{\ci}&&&&*{\ci}}
          \\[-5pt]
&
\begin{matrix}
 \\
\tilde E_6
\end{matrix}
\ \xymatrix@=10pt{
&&*{\ci}\ar@{-}[d]&&\\
&&*{\ci}\ar@{-}[d]&&\\
*{\ci}\ar@{-}[r]&*{\ci}\ar@{-}[r]&
*{\ci}
\ar@{-}[r]&*{\ci}\ar@{-}[r] &*{\ci}}
&&
\begin{matrix}
 \\
\tilde E_7
\end{matrix}
\ \xymatrix@=10pt{
&&&*{\ci}\ar@{-}[d]&&\\
*{\ci}\ar@{-}[r]&*{\ci}\ar@{-}[r]&*{\ci}\ar@{-}[r]&
*{\ci}
\ar@{-}[r]&*{\ci}\ar@{-}[r] &*{\ci}\ar@{-}[r] &*{\ci}}
          \\
&\begin{matrix}
 \\
\tilde E_8
\end{matrix}
\ \xymatrix@=10pt{
&&*{\ci}\ar@{-}[d]&&\\
*{\ci}\ar@{-}[r]&*{\ci}\ar@{-}[r]&
*{\ci}
\ar@{-}[r]&*{\ci}\ar@{-}[r] &*{\ci}\ar@{-}[r] &*{\ci}\ar@{-}[r] &*{\ci}\ar@{-}[r] &*{\ci}}
\end{aligned}
\end{equation*}

(the index plus one is the
number of vertices).

  \item Let $Q$ be of tame type
      and let ${\bf
      z}=(z_1,\dots,z_t)$ be an
      integer vector with
      nonnegative components.
      There exists an
      indecomposable
      representation of
      dimension $\bf z$ if and
      only if $q_Q({\bf z})=0$
      or $1$. (Representations
      of quivers of tame type
      were classified
      independently in
      \cite{don1} and
      \cite{naz}.)
\end{itemize}

\item 

 \emph{Kac's theorem}
\cite{kac,kraft}: Let $\mathbb F$
be an algebraically closed field.
The set of dimensions of
indecomposable representations of a
quiver $Q$ with $t$ vertices over
$\mathbb F$ coincides with the
positive root system $\Delta_+(Q)$
defined in \cite{kac}. The
following holds for ${\bf z}\in
\Delta_+(Q)$:
\begin{itemize}
  \item[(a)] $q_Q({\bf z})\le
      1$.

  \item[(b)] If $q_Q({\bf
      z})=1$, then all
      representations of
      dimension ${\bf z}$ are
      isomorphic.

  \item[(c)] If $q_Q({\bf
      z})\le 0$, then there are
      infinitely many
      nonisomorphic
      representations of
      dimension ${\bf z}$ and
      the number of parameters
      of the set of
      indecomposable
      representations of
      dimension ${\bf z}$ is
      \[1-q_Q({\bf z})=
      \sum_{i\longrightarrow
      j}z_iz_j-(z_1^2+\dots+z_t^2-1)
\text{\quad(see Example 3).}
      \]

\end{itemize}

 \item 
 {\it Belitskii's algorithm}
\cite{bel,bel1}: Let $\mathbb F$ be
an algebraically closed field.
Belitskii constructed an algorithm
that transforms each pair $(A,B)$
of $n\times n$ matrices over
$\mathbb F$ to a pair
$(A_{\text{can}},B_{\text{can}})$
that is similar to $(A,B)$ and is
such that \[ \text{$(A,B)$ is
similar to $(C,D)$}
\quad\Longleftrightarrow\quad
(A_{\text{can}},B_{\text{can}})=
(C_{\text{can}},D_{\text{can}}).\]
The pair
$(A_{\text{can}},B_{\text{can}})$
is called {\it Belitskii's
canonical form} of $(A,B)$ under
similarity. We can define {\it
Belitskii's canonical pairs} as
those matrix pairs that are not
changed by Belitskii's algorithm,
but we cannot expect to obtain an
explicit description of them.
Friedland \cite{fri} gave an
alternative approach to the problem
of classifying matrix pairs up to
similarity.

 \item 
 {\it The tame and wild theorem}
\cite{ser}: Belitskii's algorithm
was extended to a wide class of
matrix problems that includes the
problems of classifying
representations of quivers and
representations of finite
dimensional algebras. For each
matrix problem from this class over
an algebraically closed field
$\mathbb F$, denote by Bel$_{mn}$
the set of $m\times n$
indecomposable Belitskii canonical
matrices and consider Bel$_{mn}$ as
a subset in the affine space of
$m\times n$ matrices $\mathbb
F^{m\times n}$. Then
\begin{itemize}
  \item either Bel$_{mn}$
      consists of a finite
      number of points and
      straight lines for every
      $m\times n$ (then the
      matrix problem is of tame
      type),

\item or Bel$_{mn}$ contains a
    2-dimensional plane for a
    certain $m\times n$ (then
    the matrix problem is of
    wild type).
\end{itemize}
This statement is a geometric form
of Drozd's tame and wild theorem
\cite{dro}.

 \end{enumerate}
 \end{fact}\medskip

\begin{example}

 \begin{enumerate}

 \item 
Two pairs $(A,B)$ and $(A',B')$ of
$n\times n$ matrices are similar if
and only if the pairs
\[
\left(\left[ \begin{array}{lcc|r}
0&I_n&0&0\\0&0&I_n&0\\0&0&0&0\\
\hline 0&0&0&0
\end{array}  \right],\
\left[ \begin{array}{lcc|c}
0&A&0&I_n\\0&0&A&0\\0&0&0&0\\
\hline 0&0&B&0
\end{array}  \right]\right)\]
and \[ \left(\left[
\begin{array}{lcc|r}
0&I_n&0&0\\0&0&I_n&0\\0&0&0&0\\
\hline 0&0&0&0
\end{array}  \right],\
\left[ \begin{array}{lcc|c}
0&A' & 0 & I_n
\\0&0&A' &0\\0& 0 &0&0\\
\hline 0&0&B'&0
\end{array}  \right]\right)
\]
of commuting nilpotent $4n\times
4n$ matrices are similar. Thus, a
solution of the problem of
classifying pairs of commuting
nilpotent matrices up to similarity
would imply a solution of the
problem of classifying pairs of
arbitrary matrices up to
similarity.

 \item 
Two representations
\[
{\xymatrix@=15pt{
 &{\mathbb F^q}&\\
 {\mathbb F^p}
\save !<-2mm,0cm>
\ar@(ul,dl)@{->}_{A}
\restore
 \ar@{->}[ur]^{B}
 \ar@<0.4ex>[rr]^{C}
 \ar@<-0.4ex>[rr]_{D} &&{\mathbb F^r}
 \ar@{<-}[ul]_{E}
 \save
!<2mm,0cm>
\ar@{<-}@(ur,dr)^{F}
\restore }}
           \quad
\begin{matrix}\\\text{ and }\end{matrix}
           \quad
{\xymatrix@=15pt{
 &{\mathbb F^{q'}}&\\
 {\mathbb F^{p'}}
\save !<-2mm,0cm>
\ar@(ul,dl)@{->}_{A'}
\restore
 \ar@{->}[ur]^{B'}
 \ar@<0.4ex>[rr]^{C'}
 \ar@<-0.4ex>@{->}[rr]_{D'} &&{\mathbb F^{r'}}
 \ar@{<-}[ul]_{E'}
 \save
!<2mm,0cm>
\ar@{<-}@(ur,dr)^{F'}
\restore }}
\]
are isomorphic over a field
$\mathbb F$ with at least 4
distinct elements $\alpha
,\beta,\gamma ,\delta $ if and only
if the pairs
\[
\left(\begin{bmatrix}
\alpha I_{p}&0&0&0\\ 0&\beta I_{q}&0&0\\
0&0&\gamma I_{r}&0\\ 0&0&0&\delta I_{r}
\end{bmatrix},\
\begin{bmatrix}
A&0&0&0\\ B&0&0&0\\
C&0&0&0  \\ D\ &\
E\ &\ I_{r}\ &\ F
\end{bmatrix}\right) \] and \[\left(\begin{bmatrix}
\alpha I_{p'}&0&0&0\\ 0&\beta I_{q'}&0&0\\
0&0&\gamma I_{r'}&0\\ 0&0&0&\delta I_{r'}
\end{bmatrix},\
\begin{bmatrix}
A'& 0 & 0 & 0\\ B'&0&0&0\\
C'&0&0&0  \\ D'\ &\
E'\ &\ I_{r'}\ &\ F'
\end{bmatrix}\right)
\]
are similar. This example can be
extended to representations of any
quiver over any field as in
\cite[Section 5]{debora}.

 \item 

The statement about the number of
parameters in Fact 5(c) is
intuitively clear: Let $\mathcal A$
be a representation of dimension
${\bf z}$. In some bases of
$A_1,\dots,A_t$, let $M_{\alpha}$
be the matrix of $A_{\alpha}$ for
an arrow $\alpha :i\longrightarrow
j$. Let $S_1,\dots,S_t$ be the
change of basis matrices. Then the
$\sum_{i\to j}z_iz_j$ entries of
$M_{\alpha}$'s  are reduced by
$z_1^2+\dots+z_t^2$ entries of
$S_i$'s. But really only
$z_1^2+\dots+z_t^2-1$ independent
parameters are used since
multiplying  all $S_i$ by the same
nonzero scalar does not change the
transformation $M_{\alpha}\mapsto
S_j^{-1}M_{\alpha}S_i$  for all
arrows $\alpha :i\longrightarrow
j$.

\end{enumerate}
\end{example}

\vspace*{1pc}
\section{Quivers of finite dimensional
algebras}\label{ser3}

All representations of a finite
dimensional algebra can be identified
with all representations of some quiver
with relations.

\begin{definition}

\noindent A {\bf relation} in a quiver
$Q$ over a field $\mathbb F$ is a
formal expression of the form
\begin{equation}\label{a1a}
\sum_{i=1}^m
c_i\alpha_{ip_i}\cdots
\alpha_{i2}\alpha_{i1}
=0,\qquad 0\ne c_i\in \mathbb F,
\end{equation}
in which all
\[ u\xrightarrow[\phantom{\qquad}]
 {\alpha_{i1}}
  u_{i2}
  \xrightarrow[\phantom{\qquad}]
  {\alpha_{i2}}
         \cdots
  \xrightarrow[\phantom{\qquad}]
  {\alpha_{i,p_i-1}}
  u_{ip_i} \xrightarrow[\phantom{\qquad}]
  {\alpha_{ip_i}}
  v,\qquad i=1,\dots,m,
\]
are {\bf directed paths} in $Q$ with
the same start vertex $u$ and the same
end vertex $v$ (it is possible that
$u=v$).

A representation $\mathcal A$ of $Q$
{\bf satisfies the relation
\eqref{a1a}} if
\begin{equation}\label{1mau}
\sum_{i=1}^m c_iA_{\alpha_{ip_i}}\cdots
A_{\alpha_{i2}}
A_{\alpha_{i1}}
=0.
\end{equation}
If $u=v$, then \eqref{a1a} may have a
summand $c_i\varepsilon _u$, in which
$\varepsilon _u$ is the path without
arrows. This ``lazy'' path $\varepsilon
_u$ (to stand in place) is replaced in
\eqref{1mau} by the identity operator
on $A_u$.

By a {\bf quiver with relations}
$(Q,L)$ we mean a quiver $Q$ with a
finite set $L$ of relations in $Q$. Its
{\bf set of representations} consists
of all representations of $Q$ that
satisfy all relations from $L$.

The {\bf path algebra} $\mathbb FQ$ of
a quiver $Q$ is a finite dimensional
algebra over a field $\mathbb F$ whose
elements are formal linear combinations
\[\sum_{i=1}^m c_i\alpha_{ip_i}\cdots
\alpha_{i2}\alpha_{i1},\] in which
$c_i\in \mathbb F $ and
$\alpha_{ip_i}\cdots
\alpha_{i2}\alpha_{i1}$ are directed
paths (they may be lazy paths and may
have distinct start vertices and
distinct end vertices). Their
multiplication is determined by the
distributive law and the rule:
\begin{multline*}
(\beta_q\cdots
\beta_{1})(\alpha_p\cdots
\alpha_{1})\\=
\begin{cases}
\beta_q\cdots
\beta_{1}\alpha_p\cdots
\alpha_{1}& \hbox{if the end vertex of $\alpha_p$ is the start vertex of $\beta_{1}$}, \\
               0 & \hbox{otherwise}.
             \end{cases}
\end{multline*}
The multiplicative identity of the
algebra $\mathbb FQ$ is the sum
$\varepsilon _1+\dots+\varepsilon _t$
of all lazy paths taken over all
vertices. If $(Q,L)$ is a quiver with
relations, then its {\bf path algebra}
$\mathbb F(Q,L)$ is determined modulo
these relations; that is, $\mathbb
F(Q,L):=\mathbb FQ/{\mathcal L}$ in
which ${\mathcal L}$ is the two-sided
ideal of $\mathbb FQ$ generated by the
left-hand sides of relations from $L$.

A {\bf representation} of a finite
dimensional algebra $\Lambda$ over
$\mathbb F$ is a homomorphism $\varphi:
\Lambda\to \End V$ to the algebra $\End
V$ of linear operators on a vector
space $V$ over $\mathbb F$.

An algebra over an algebraically closed
field $\mathbb F$ is called  a {\bf
basic algebra} if for some positive
integer $m$ it is isomorphic to an
algebra $\Lambda $ of upper triangular
$m\times m$ matrices over $\mathbb F$
that satisfies the condition:
\[
\begin{bmatrix}
  a_{11}&\dots & a_{1m} \\
   & \ddots&\vdots \\0&&a_{mm}
\end{bmatrix}\in\Lambda \qquad
\Longrightarrow
\qquad\begin{bmatrix}
   a_{11} & & 0 \\
   & \ddots& \\0&& a_{mm}
\end{bmatrix}\in\Lambda .
\]
\end{definition}
\medskip

\begin{fact}

 \begin{enumerate}
 \item 
 Each finite dimensional algebra
 $\Lambda $ over a field $\mathbb
 F$ is isomorphic to the path
 algebra $\mathbb F(Q,L)$ of a
 quiver with relations $(Q,L)$,
 which is
constructed in
 \cite[Chapter~II]{simson}. We give
 a simplified construction of
 $(Q,L)$ in the
 following algorithm.\\

\noindent
\fbox{\parbox{12.7cm}{{\bf
Algorithm 1: From an algebra
\pmb{$\Lambda $} to a quiver
\pmb{$(Q,L)$}.}

 \begin{itemize}
 \item[1.] Decompose the unit
     of $\Lambda $ into a sum
     of orthogonal idempotents:
\begin{equation}\label{jyk}
1=e_1+\dots+e_t,\qquad  e_ie_j=0\text{ if }i\ne j,\ \
e_i^2=e_i\ne 0.
\end{equation}

\item[2.] Choose a set
    $a_1,\dots,a_n$ of elements
    of $\Lambda $ such that
    $e_1,\dots,e_t,a_1,\dots,a_n$
    generate $\Lambda $ and
    each $a_i$ is equal to
    $e_{q(i)}a_ie_{p(i)}$ for
    some $p(i)$ and $q(i)$
    (such a set exists since if
    $b_1,b_2,\dots$ generate
    $\Lambda $, then all
    $e_ib_je_k$ also generate
    $\Lambda $).

\item[3.] Denote by $Q$ the
    quiver with vertices
    $1,\dots,t$ and $n$ arrows
    $\alpha
    _i:p(i)\longrightarrow
    q(i)$.

\item[4.] Denote by $\pi :
    \mathbb FQ\to \Lambda$ the
    epimorphism of algebras
    such that
    \[\pi(\varepsilon_1)=e_1,\ \dots,\
\pi(\varepsilon_t)=e_t,\ \
    \pi(\alpha_1)=a_1,\ \dots,\
    \pi(\alpha_n)=a_n.\]

\item[5.] Construct a set $L$
    of relations in $Q$ by
    choosing a finite subset of
    $\cup \varepsilon
    _j\Ker(\pi)\varepsilon _i$
    that generates $\Ker(
    \pi)$, expressing its
    elements through
    $\varepsilon
    _1,\dots,\varepsilon
    _t,\alpha _1,\dots,\alpha
    _n$, and equating them to
    zero.
\end{itemize}
Then $\mathbb F(Q,L)\simeq \mathbb
FQ/\Ker(\pi)\simeq \Lambda $.
}}\bigskip

\item 

In the following algorithm, we
construct a canonical
correspondence
$\varphi\longleftrightarrow \cal R$
between representations of
$\Lambda$ and representations of
the quiver $(Q,L)$ constructed by
Algorithm~1 such that $\varphi$ and
$\varphi'$ are isomorphic if and
only if the corresponding
representations $\cal R$
and ${\cal R}'$ are isomorphic.\\

\noindent
\fbox{\parbox{12.7cm}{{\bf
Algorithm 2: From a representation
$\mathbf{\varphi: \Lambda\to
\text{End}\,V}$ of \pmb{$\Lambda$}
to a representation \pmb{$\cal R$}
of \pmb{$(Q,L)$}.}
 \begin{enumerate}
 \item[1] Since \eqref{jyk}
     holds for $\tau_i:=\varphi
     (e_i):V\to V$ instead of
     $e_i$, we have
     $V=\tau_1V\oplus\dots\oplus
     \tau_tV$. Put
     $R_i:=\tau_iV$ for every
     vertex $i=1,\dots,t$.

\item[2] For each $a_i$ from
    Step 2 of Algorithm~1,
    define $\rho_i:=\varphi
    (a_i):V\to V$. Since
    $\rho_i=\tau_{q(i)}
    \rho_i\tau_{p(i)}$, we have
    $\rho_i(\tau
    _{p(i)}V)\subseteq
    \tau_{q(i)}V$ and
    $\rho_i(\tau_kV)=0$ if
    $k\ne p(i)$, and so each
    $\rho_i$ is fully
    determined by its
    restriction $\rho_i|_{\tau
    _{p(i)}V}: \tau _{p(i)}V\to
    \tau_{q(i)}V$. Put
    $R_{\alpha_i}:=\rho
    _i|_{\tau _{p(i)}V}$ for
    every arrow $\alpha
    _i:p(i)\longrightarrow
    q(i)$.
\end{enumerate}
}}
\bigskip

\item 
If the field $\mathbb F$ is
algebraically closed, then it
suffices to study representations
of basic algebras since for each
finite dimensional algebra over
$\mathbb F$ there exists a basic
algebra over $\mathbb F$ such that
the categories of representations
of these algebras are equivalent;
see \cite[Corollary
I.6.10]{simson}. One usually
applies Algorithms 1 and 2 to a
basic algebra $\Lambda$ over
$\mathbb F$, chooses
$a_1,\dots,a_n$ among its nilpotent
elements, and takes the numbers $t$
and $n$ to be maximal and minimal,
respectively.
\end{enumerate}

 \end{fact}

\begin{example}

 \begin{enumerate}

 \item 
The path algebra $\mathbb F(Q,L)$
of the quiver with relation
\[\raisebox{14pt}{\xymatrix@=0pt{
&&&{2}\ar@{<-}[llld]_{\alpha }
\ar[rrrd]^{\beta }\\
{1}&&&&&&{4}
\\&&&{3}\ar@{<-}[lllu]^{\gamma }\ar[rrru]_{\delta
}}}\qquad\qquad \beta \alpha=\delta
\gamma\] has the basis
\[\varepsilon_1,\ \varepsilon_2,\
\varepsilon_3,\ \varepsilon_4,\
\alpha,\ \beta,\ \gamma,\ \delta,\ \beta
\alpha \]  over $\mathbb F$. The
product of $\varepsilon _3-\delta
+\beta \alpha$ and $\varepsilon
_1+2\gamma$ in $\mathbb F(Q,L)$ is
\[\varepsilon _3\varepsilon
_1-\delta \varepsilon _1+\beta
\alpha\varepsilon _1+2\varepsilon
_3\gamma-2\delta\gamma +2\beta
\alpha\gamma= \beta\alpha+2\gamma
-2\delta \gamma =-\beta \alpha
+2\gamma.\] Each representation
$\mathcal A$ of $(Q,L)$ defines a
representation of $\mathbb F(Q,L)$
by operators on the space
$A_1\oplus A_2\oplus A_3\oplus
A_4$.

 \item 
The problem of classifying
representations of the quiver with
relations
\begin{equation*}\label{grj}
\xymatrix{
 {1} \ar@(ul,dl)@{->}_{\alpha}
 \ar@(ur,dr)@{<-}^{\beta}}\qquad\qquad
 \alpha \beta =\beta \alpha =0
\end{equation*}
is the problem of classifying pairs
of mutually annihilating linear
operators, which was solved in
\cite{gel-pon1} (see also
\cite{BGS09}). Its path algebra is
an infinite dimensional algebra
whose elements are finite linear
combinations of $\varepsilon
_1,\alpha, \alpha^2, \alpha^3,
\dots,\beta, \beta ^2, \beta ^3,
\dots$ over $\mathbb F$.

\item Let us apply Algorithm 1 to
    the basic algebra
\[
\Lambda:=\Bigl\{\ \begin{bmatrix}
\;u\;&\;x\;&\;z\;\\0&u&y\\0&0&v
\end{bmatrix}\,:\,u,v,x,y,z\in\mathbb C
\ \Bigr\}.
\]
Decompose its unit into a sum of
      orthogonal idempotents:
$I_3=e_1+e_2$, in which
$$e_1:=\diag(1,1,0),\qquad
e_2:=\diag(0,0,1).$$ Write
\[
a_1:=\begin{bmatrix}
\;0\;&\;1\;&\;0\;\\0&0&0\\0&0&0
\end{bmatrix}=e_1a_1e_1,\qquad
a_2:=\begin{bmatrix}
\;0\;&\;0\;&\;0\;\\0&0&1\\0&0&0
\end{bmatrix}=e_1a_2e_2\,.
\]
The elements $e_1,e_2,a_1,a_2$
generate $\Lambda $. We obtain the
quiver with relations
\begin{equation}\label{kdt}
\xymatrix{
 {1} \ar@(ul,dl)@{->}_{\alpha_1}
 &2\ar[l]_{\alpha_2}}\qquad\qquad
 \alpha_1^2=0
\end{equation}

whose path algebra is isomorphic to
$\Lambda $. Each representation of
$\Lambda $ is obtained from a
representation of \eqref{kdt} and
vice versa.
\end{enumerate}
\end{example}

\vspace*{1pc}
\section{Systems of linear mappings and
forms as representations of mixed
graphs}\label{ser4}

By analogy with quiver representations,
systems of linear mappings and forms
can be considered as representations of
mixed graphs, in which forms are
assigned to undirected edges.

\begin{definition}

\noindent Let $\mathbb F$ be a field
with a fixed involution $a\mapsto
\overline a$; i.e., a bijection
$\mathbb F\to \mathbb F$ (which can be
the identity) satisfying
\[\overline{a + b} = \overline a +
\overline b,\qquad \overline{ab} =
\overline a\overline b,\qquad
\overline{\overline a}=a.\]

A {\bf mixed graph} $G$ is a graph in
which loops and multiple edges are
allowed and that may contain both
directed and undirected edges; we
suppose that the vertices are
$1,\dots,t$.

A {\bf representation} $\mathcal A$ of
$G$ over $\mathbb F$ is given by
assigning to each vertex $i$ a finite
dimensional vector space $A_i$ over
$\mathbb F$, to each directed edge
$\alpha:i\longrightarrow j$ a linear
mapping $A_{\alpha}:A_i\to A_j$, and to
each undirected edge $\lambda :i\lin j$
($i\le j$) a sesquilinear form
$A_{\lambda }: A_j\times A_i\to \mathbb
F$; this form is bilinear if the fixed
involution on $\mathbb F$ is the
identity. We suppose that $A_{\lambda
}$ is linear on $A_i$ and semilinear on
$A_j$.

An {\bf isomorphism} $\varphi: \mathcal
A \ \begin{matrix}\sim\; \\[-9pt] \to
\\[-9pt]{}\end{matrix}\; \mathcal
B$ between representations $\mathcal A$
and $\mathcal B$ of $G$ is a family of
linear bijections \[\varphi_1:A_1\to
B_1,\ \dots,\ \varphi_t:A_t\to B_t\]
such that $$\text{$\varphi_jA_{\alpha
}=B_{\alpha }\varphi _i$\qquad for each
directed edge $\alpha :i\longrightarrow
j$}$$ and $$\text{$A_{\lambda
}(y,x)=B_{\lambda
}(\varphi_jy,\varphi_ix)$\qquad for
each undirected edge $\lambda :i\lin
j$.}$$

The notions of the {\bf dimension} of a
representation, {\bf direct sum}, and
{\bf finite, tame, and wild types} are
defined for mixed graphs in the same
way as for quivers. The {\bf Tits form}
is defined as in \eqref{kur}, but the
sum is also taken over all undirected
edges $i\lin j$ ($i\le j$).

\end{definition}

\begin{fact}

 \begin{enumerate}

 \item 
The Krull--Schmidt, Gabriel, and
Donovan--Freislich--Nazarova
theorems (the Fact in Section
\ref{ser1} and Facts 3 and 4 in
Section \ref{ser2}) remain true if
we replace the word ``quiver'' by
``mixed graph''; see
\cite{ser_izv}.

 \item 
 {\it A generalization of Kac's
theorem} (from Fact 5 of Section
\ref{ser2})
\cite{ser1dir,ser_brazil}: Let
$\mathbb F$ be an algebraically
closed field of characteristic not
$2$. The set of dimensions of
indecomposable representations of a
mixed graph $G$ over $\mathbb F$
coincides with the positive root
system $\Delta_+(G)$ (its
definition in \cite{kac} does not
depend on the orientation of
edges).
 \end{enumerate}
 \end{fact}

\begin{exampl}

 \begin{enumerate}
 \item[]  Each representation
     \[\raisebox{20pt}
     
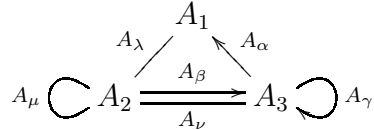
{\xymatrix@R=15pt@C=10pt{
 &{A_1}&\\
 {A_2}
\save !<-2mm,0cm>
\ar@(ul,dl)@{-}_{A_{\mu}}
\restore
 \ar@{-}[ur]^{A_{\lambda}}
\ar@<0.4ex>[rr]^{A_{\beta}}
 \ar@<-0.4ex>@{-}[rr]_{A_{\nu}}
 &&{A_3} \ar[ul]_{A_{\alpha}} \save
 !<2mm,0cm>
\ar@(ur,dr)^{A_{\gamma}} \restore
}}\] {of}\[
 \raisebox{23pt}
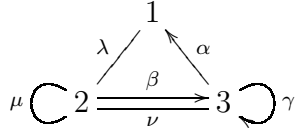
{\xymatrix@R=20pt@C=15pt{
 &{1}&\\
 {2}\ar@(ul,dl)@{-}_{\mu}
 \ar@{-}[ur]^{\lambda}
\ar@<0.4ex>[rr]^{\beta}
 \ar@<-0.4ex>@{-}[rr]_{\nu} &&{3}
 \ar[ul]_{\alpha}
 \ar@(ur,dr)^{\gamma} }} \] is a
 system of vector spaces
$A_1,A_2,A_3$ over $\mathbb F$,
linear mappings $A_{\alpha}$,
$A_{\beta}$, $A_{\gamma}$, and
sesquilinear forms \[A_{\lambda}:
A_2\times A_1\to {\mathbb F},\quad
A_{\mu}: A_2\times A_2\to {\mathbb
F},\quad A_{\nu}: A_3\times A_2\to
 {\mathbb F}\] (these forms are
 bilinear if the fixed involution
 on $\mathbb F$
is the identity).
\end{enumerate}
\end{exampl}

\vspace*{1pc}
\section{Generalization of the law
of inertia to representations of mixed
graphs}\label{ser5}

The problem of classifying systems of
forms and linear mappings over $\mathbb
C$ and $\mathbb R$ is reduced to the
problem of classifying systems of
linear mappings.

\begin{definition}

\noindent Let $V$ be a finite
dimensional vector space over a field
$\mathbb F$ with a fixed involution
(which can be the identity). By the
{\bf $^{\filledstar\!}$dual space} of
$V$, we mean the space
$V^{\filledstar}$ of all mappings
$\varphi: V\to \mathbb F$ that are {\bf
semilinear}, i.e., \[
\varphi(au+bv)=\bar{a}\varphi (u)+
\bar{b}\varphi (v)\] for all $u,v\in V$
and $a,b\in \mathbb F$.

For each linear mapping $A:U\to V$, we
define the {\bf
$^{\filledstar\!}$adjoint mapping}
$$A^{{\filledstar}}:
V^{{\filledstar}}\to U^{{\filledstar}}$$
by putting
$$A^{{\filledstar}}\varphi:=\varphi
A\qquad\text{
for all $\varphi\in V^{\filledstar}$.}$$

For each mixed graph $G$, we denote by
$\underline{G}$ the quiver that is
obtained from $G$ by replacing
\begin{itemize}\parskip=-3pt
  \item each vertex $i$ of ${G}$ by
      the vertices $i$ and
      $i^{\filledstar}$,
  \item each arrow $\alpha:
      i\longrightarrow j$ by the
      arrows \[\alpha:
      i\longrightarrow j,\qquad
      \alpha^{\filledstar}:
      j^{\filledstar}\longrightarrow
      i^{\filledstar},\]
  \item each undirected edge
      $\lambda : i\lin\, j\ (i\le
      j)$ by the arrows \[\lambda :
      i\longrightarrow
      j^{\filledstar},\qquad
      \lambda^{\filledstar}:
      j\longrightarrow
      i^{\filledstar}.\]
\end{itemize}
We consider $\underline{G}$ as a {\bf
quiver with involution} on the set of
vertices and on the set of arrows.

For each representation $\mathcal A$
over $\mathbb F$ of a mixed graph $G$,
we denote by $\underline {\mathcal A}$
the representation of $\underline {G}$
that is obtained from $\mathcal A$ by
replacing (see Example \ref{jjy})
\begin{itemize}
 \parskip=-3pt
  \item each vector space $A_i$ by
      the mutually
      $^{\filledstar\!}$dual spaces
      $A_i$ and
      $A_i^{\filledstar}$,
  \item each linear mapping
      $A_{\alpha}: A_i\to A_j$ by
      the mutually
      $^{\filledstar\!}$adjoint
      mappings \[A_{\alpha}: A_i\to
      A_j,\qquad
      A_{\alpha}^{\filledstar}:
      A_j^{\filledstar}\to
      A_i^{\filledstar},\]

  \item each sesquilinear form
      $A_{\lambda}: A_j\times
      A_i\to \mathbb F$ by the
      mutually
      $^{\filledstar\!}$adjoint
      mappings \[ A_{\lambda}: u\in
      A_i\mapsto
      A_{\lambda}(?,u)\in
      A_j^{\filledstar},\qquad
      A_{\lambda}^{\filledstar}:
      v\in A_j\mapsto
      \overline{A_{\lambda}(v,?)}\in
      A_i^{\filledstar}.\]
\end{itemize}

For each representation $\mathcal M$ of
$\underline {G}$, we define the {\bf
adjoint representation} ${\mathcal
M}^{\circ}$ of $\underline {G}$ that is
formed by the vector spaces
${M}^{\circ}_v:={M}^{\filledstar}_{v^{\filledstar}}$
for all vertices $v$ of $\underline{G}$
and the linear mappings
${M}^{\circ}_{\tau}:={M}^{\filledstar}_{\tau^{\filledstar}}$
for all arrows $\tau$ of
$\underline{G}$ (see Example
\ref{jir}).

A representation $\mathcal M$ of
$\underline {G}$ is {\bf selfadjoint}
if $\mathcal M^{\circ} =\mathcal M$.

A {\bf mixed graph with relations}
$(G,L)$ is a mixed graph $G$ with a
finite set $L$ of relations in
$\underline G$. By {\bf representations
of} $(G,L)$ we mean those
representations $\mathcal A$ of $G$ for
which $\underline{\mathcal A}$
satisfies $L$.

For each relation \[\sum_{i=1}^m
c_i\tau_{ip_i}\cdots
\tau_{i2}\tau_{i1}=0\qquad\text{in
$\underline {G}$ (see \eqref{a1a}),}\]
we define the {\bf adjoint relation}
\[\sum_{i=1}^m
\bar{c}_i\tau_{i1}^{\filledstar}
\tau_{i2}^{\filledstar}\cdots
\tau_{ip_i}^{\filledstar}
=0\qquad\text{in $\underline {G}$}.\]
For each set $L$ of relations in
$\underline G$, we denote by
$L^{\filledstar}$ the set of relations
that are adjoint to the relations from
$L$.

For each representation ${\mathcal
 A}$ of $G$,
we denote by ${\mathcal A}^-$ the
representation of $G$ that is obtained
from ${ \mathcal A}$ by replacing all
its forms ${ A}_{\lambda }$ by $-{
A}_{\lambda }$.

\end{definition}

\begin{fac}
\begin{enumerate}
\item[]
\noindent In the following
algorithm, the problem of
classifying representations of a
mixed graph with relations $(G,L)$
over $\mathbb C$ and $\mathbb R$ is
reduced to the problem of
classifying representations of the
quiver with relations
$(\underline{G}, L\cup
L^{\filledstar})$. The algorithm is
a special case of the method
\cite{ser_izv} (see also
\cite{roi,hor+ser_anyfield,
ser_brazil,ser_isom}) for reducing
the problem of classifying
representations of a mixed graph
$(G,L)$ over a field or skew field
$\mathbb F$ of characteristics not
2 to the problem of classifying
representations of the quiver
$(\underline{G}, L\cup
L^{\filledstar})$ over $\mathbb F$
and the problem of classifying
Hermitian and symmetric forms over
fields and skew fields that
are finite extensions of the center of $\mathbb F$.\\

\noindent
\fbox{\parbox{12.7cm}{{\bf
Algorithm 3: Classification of
representations of a mixed graph
with relations \pmb{$(G,L)$}.}

 \begin{itemize}
 \item[1.] %
Construct a set
$\ind(\underline{G}, L\cup
L^{\filledstar})$ of
indecomposable representations
of $(\underline{G}, L\cup
L^{\filledstar})$ such that
every indecomposable
representation of
$(\underline{G}, L\cup
L^{\filledstar})$ is isomorphic
to exactly one representation
from $\ind (\underline{G},
L\cup L^{\filledstar})$.

\item[2.] %
Improve $\ind (\underline{G},
L\cup L^{\filledstar})$ such
that
\begin{itemize}\parskip=-2pt
  \item[$\bullet$] if
      $\mathcal M\in\ind
      (\underline{G}, L\cup
      L^{\filledstar})$ is
      isomorphic to a
      selfadjoint
      representation, then
      $\mathcal M$ is
      selfadjoint,

  \item[$\bullet$] if
      $\mathcal M\in\ind
      (\underline{G}, L\cup
      L^{\filledstar})$ is
      not isomorphic to
      $\mathcal
M^{\circ}$, then $\mathcal
M^{\circ}\in\ind
      (\underline{G}, L\cup
      L^{\filledstar})$.
\end{itemize}
\end{itemize}

Then every representation of $(G,
L)$ over $\mathbb F$ is isomorphic
to a direct sum, uniquely
determined up to permutation of
summands, of representations of the
types
\[\begin{array}{rll}
\mathcal
A,\
\mathcal B\ &
\text{if $\mathbb F=\mathbb C$ with the indentity involution,}
                     \\
\mathcal A,\
\mathcal B,\ \mathcal B^- & \text{if $\mathbb F=\mathbb C$ with complex conjugation,}
                      \\
\mathcal A,\
    \mathcal B, \text{ and also } {\mathcal B}^{-}
\text{ if }\mathcal B^-\not\simeq \mathcal B \
  &
\text{if $ \mathbb F=\mathbb R$},
\end{array}
\]
in which $\underline {\mathcal
A}=\mathcal M\oplus \mathcal
M^{\circ}$ for each unordered pair
$\{\mathcal M,\mathcal M^{\circ}\}$
such that $\mathcal
M^{\circ}\ne\mathcal M\in\ind
(\underline{G}, L\cup
L^{\filledstar})$ and $\underline
{\mathcal B}\in\ind (\underline{G},
L\cup
L^{\filledstar})$. }}\\

Thus, each system of linear
mappings and bilinear forms over
$\mathbb C$ or $\mathbb R$ and each
system of linear mappings and
sesquilinear forms over $\mathbb C$
are decomposed into direct sums of
indecomposables uniquely, up to
isomorphisms of summands. This is
the Krull--Schmidt theorem for
representations of mixed graphs;
see Fact 1 in Section \ref{ser4}.

\end{enumerate}
 \end{fac}

\begin{example}
 \begin{enumerate}
 \item 
The problems of classifying
representations over a field
$\mathbb F$ of the mixed graphs
with relations
\begin{itemize}
  \item
      $\,\phantom{{\scriptstyle
      \lambda}\mrof} 1\form
      {\scriptstyle\lambda }$

  \item $\,{\scriptstyle
      \lambda}\mrof 1\form
      {\scriptstyle\mu}\quad
      \lambda =\varepsilon
      \lambda ^{\filledstar},\
      \ \mu =\delta \mu
      ^{\filledstar}$

  \item $ {\,\scriptstyle
      \alpha
      }\!\lefttorightarrow \!\!
      1\form
      {\scriptstyle\lambda
      }\quad \lambda
 =\varepsilon \lambda
 ^{\filledstar}\text{ is
 nonsingular, }\ \
 \alpha^{\filledstar}\lambda=\lambda
 \alpha$

  \item ${\,\scriptstyle \alpha
      }\!\lefttorightarrow \!\!
      1\form
      {\scriptstyle\lambda
      }\quad \lambda
      =\varepsilon \lambda
      ^{\filledstar}
      =\alpha^{\filledstar}\lambda
      \alpha\text{ is
      nonsingular}$
\end{itemize}
(in which
$\varepsilon,\delta\in\{-1,1\}$)
are the problems of classifying
\begin{itemize}\parskip=-2pt
  \item sesquilinear (bilinear
      if the involution on
      $\mathbb F$ is the
      identity) forms,

  \item pairs of
      $\varepsilon$-,
      $\delta$-Hermitian
      (symmetric, or
      skew-symmetric)  forms,

  \item triples $(V,H,A)$, in
      which $V$ is a vector
      space, $H$ is a
      nonsingular
      $\varepsilon$-Hermitian
      (symmetric, or
      skew-symmetric) form on
      $V$, and $A$ is a linear
      operator on $V$ that is
      \emph{$^H\!\!$selfadjoint},
      i.e.,
      $$H(Ax,y)=H(x,Ay)\qquad\text{for
      all $x,y\in V$}.$$

  \item triples $(V,H,A)$, in
      which $V$ is a vector
      space, $H$ is a
      nonsingular
      $\varepsilon$-Hermitian
      (symmetric, or
      skew-symmetric) form on
      $V$, and $A$ is a linear
      operator on $V$ that is
      \emph{$^H\!$unitary},
      i.e.,
      $$H(Ax,Ay)=H(x,y)\qquad\text{for
      all $x,y\in V$}.$$
\end{itemize}
Canonical matrices for these
problems are given in
\cite{ser_izv} over any field
$\mathbb F$ of characteristic not 2
up to classification of Hermitian
and symmetric forms over finite
extensions of $\mathbb F$; they are
based on the elementary divisors
rational canonical form (see
Example 2 in Section \ref{ser1}).
Simpler canonical matrices over
$\mathbb C$ and $\mathbb R$ that
are based on the Jordan canonical
form are given in
\cite{tridiag,goh1a,HS06,
hor+ser_anyfield,ser_isom,thom}.

 \item 
The problem of classifying
representations of a mixed graph
$(G,L)$ is hopeless if the quiver
$(\underline G,L\cup
L^{\filledstar})$ is of wild type.
For example, the problem of
classifying triples of Hermitian
forms and the problem of
classifying normal operators on a
complex space with scalar product
given by a nonsingular Hermitian
form are hopeless.

 \item \label{jjy} If $$\mathcal
     A:\raisebox{16pt}{\xymatrix@R=4pt{
     **[r]{A_1}\ar[dd]_{A_{\alpha
     }}
 \ar@{-}@/^/[dd]^{A_{\lambda }}\\
  \\
**[r] {A_2}\form {\scriptstyle
A_{\mu }}}} $$ is a representation
of the mixed graph $$G:
\raisebox{16pt}{\xymatrix@R=4pt{
**[r] {1}\ar[dd]_{\alpha}
 \ar@{-}@/^/[dd]^{\lambda }\\
  \\
**[r]{\: 2\form {\scriptstyle \mu
}}}}, $$ then
\begin{equation}\label{4.1}
\raisebox{20pt}{\xymatrix@R=4pt{
 &{A_1}\ar[dd]_{A_{\alpha }}
 \ar[ddrr]^(.25){A_{\lambda }}&
 &{A_1^{\filledstar}} \\
 {\underline{\mathcal A}:}\!\!\!\!\!\!\!\!\\
 &{A_2}\ar[uurr]^(.75){A_{\lambda }^{\filledstar}}
 &&{A_2^{\filledstar}}\ar[uu]_{A_{\alpha }^{\filledstar}}
 \ar@{<-}@<0.4ex>[ll]^{A_{\mu }^{\filledstar}}
 \ar@{<-}@<-0.4ex>[ll]_{A_{\mu }}
 }}
\qquad
\raisebox{20pt}{\xymatrix@R=4pt{
 &{1}\ar[dd]_{{\alpha}}
 \ar[ddrr]^(.25){\lambda }&
 &{1^{\filledstar}} \\
 {\underline{G}: }\!\!\!\!\!\!\!\!\\
 &{2}\ar[uurr]^(.75){\lambda ^{\filledstar}}
 &&{2^{\filledstar}}\ar[uu]_{{\alpha}^{\filledstar}}
\ar@{<-}@<-0.4ex>[ll]_{\mu}
 \ar@{<-}@<0.4ex>[ll]^{{\mu}^{\filledstar}}
 }}
\end{equation}\\[-15pt]

\item \label{jir} For each
    representation $\cal M$ of the
    quiver ${\underline{G}}$ in
\eqref{4.1}, the adjoint
representation $\mathcal M^{\circ}$
is constructed as follows:
\[
\xymatrix@R=4pt{
 &{U_1}\ar[dd]_{A_1}
 \ar[ddrr]^(.25){B_1}&
 &{U_2} \\
 {\mathcal M:}\!\!\!\!\!\!\!\!\\
 &{V_1}\ar[uurr]^(.75){B_2}
 &&{V_2}\ar[uu]_{A_2}
 \ar@{<-}@<0.4ex>[ll]^{C_2}
 \ar@{<-}@<-0.4ex>[ll]_{C_1}
 } \qquad
 \xymatrix@R=4pt{
 &{U_2^{\filledstar}}\ar[dd]_{A_2^{\filledstar}}
 \ar[ddrr]^(.25){B_2^{\filledstar}}&
 &{U_1^{\filledstar}} \\
 {\mathcal M^{\circ}:}\!\!\!\!\!\!\!\!\\
 &{V_2^{\filledstar}}
 \ar[uurr]^(.75){B_1^{\filledstar}}
 &&{V_1^{\filledstar}}\ar[uu]_{A_1^{\filledstar}}
 \ar@{<-}@<0.4ex>[ll]^{C_1^{\filledstar}}
 \ar@{<-}@<-0.4ex>[ll]_{C_2^{\filledstar}}
 }
\]

 \item %
Applying Algorithm 3 to Hermitian
or symmetric forms gives the law of
inertia. Indeed, these forms are
representations of the mixed graph
with relations $$(G,L): 1\form
{\scriptstyle
\lambda=\lambda^{\filledstar}}.$$
Its quiver $$(\underline G,L\cup
L^{\filledstar}):
     \!\!\xymatrix@=15pt{{1}
 \ar@<0.4ex>[r]^{\scriptscriptstyle\lambda}
 \ar@<-0.4ex>
 [r]_{\scriptscriptstyle\lambda^{\filledstar}}
 &{1^{\filledstar}}}\!\!,\
 \lambda^{\filledstar}=\lambda. $$
 By \eqref{key}, $\ind(\underline
G,L\cup L^{\filledstar})$
consists of 3 representations:\\[-5pt]
\[
\begin{matrix}
\mathcal M: \\[4mm]
\end{matrix}
\ \xymatrix{ 0
\ar@<0.4ex>[r]^{0_{10}}
 \ar@<-0.4ex>[r]_{0_{10}}&
\mathbb F},\qquad
\begin{matrix}
\mathcal M^{\circ}: \\[4mm]
\end{matrix}
\
\xymatrix{
\mathbb F
\ar@<0.4ex>[r]^{0_{01}}
 \ar@<-0.4ex>[r]_{0_{01}}&
0},\qquad
\begin{matrix}
\mathcal N=\mathcal N^{\circ}: \\[4mm]
\end{matrix}
\ \xymatrix{ \mathbb F
\ar@<0.4ex>[r]^{[1]}
 \ar@<-0.4ex>[r]_{[1]}&
\mathbb F}.
\]
Since $\mathcal M\oplus \mathcal
M^{\circ}=\underline{\mathcal  A}$
for $\mathcal A: {\mathbb F}\form {
[0]}$ and $\mathcal
N=\underline{\mathcal B}$ for
$\mathcal B: {\mathbb
F}\form{[1]}$, Algorithm 3 ensures
that each representation of $
 {1}\form{\scriptstyle
 \lambda=\lambda^{\filledstar}}$ is
 isomorphic to a direct sum,
 uniquely determined up to
 permutation of summands, of
 representations of the form:
\[\begin{array}{lll}
 {\mathbb C}\form [0],\ \
 {\mathbb C}\form [1]
 & \text{if $\mathbb F=\mathbb C$
 with the identity involution,}
                     \\
 {\mathbb C}\form [0],\ \
 {\mathbb C}\form [1],\ \
 {\mathbb C}
 \form [-1]
& \text{if $\mathbb F=\mathbb C$
with complex conjugation
or $\mathbb R$.}
\end{array}
\]

 \item %
(See details in
\cite{HS06,ser_izv}.) Applying
Algorithm 3 to sesquilinear or
bilinear forms over $\mathbb C$
gives their canonical forms from
\cite{HS06}. Indeed, these forms
are representations of the graph
$G: 1\form {\scriptstyle\lambda }$.
Each representation {${\cal
M}:\xymatrix@1@=20pt{U
 \ar@<0.4ex>[r]^{\scriptstyle A\ }
 \ar@<-0.4ex> [r]_{\scriptstyle B\
 } &{\: V}}$} of the quiver
$ \underline
G:\xymatrix@1@=20pt{{1}
 \ar@<0.4ex>[r]^{\scriptscriptstyle\lambda\
 } \ar@<-0.4ex>
 [r]_{\scriptscriptstyle\lambda^{\filledstar}\
 } &{\: 1^{\filledstar}}}$  defines the
 representations\\[-5pt]
\begin{equation*}\label{hge}
\mathcal
M^{\circ}:\xymatrix@1@=20pt{V^{\filledstar}
 \ar@<0.4ex>[r]^{\scriptstyle
 B^{\filledstar}} \ar@<-0.4ex>
 [r]_{\scriptstyle A^{\filledstar}} &{\:
 U^{\filledstar}}},\qquad
 {\mathcal
 M}\oplus\mathcal
M^{\circ}:\ \xymatrix@1@=40pt{U\oplus
V^{\filledstar}
 \ar@<0.4ex>[r]^{
\left[\begin{smallmatrix}
  0 & B^{\filledstar} \\
  A & 0 \\
\end{smallmatrix}\right]\ } \ar@<-0.4ex>
 [r]_{
\left[\begin{smallmatrix}
  0 & A^{\filledstar} \\
  B & 0 \\
\end{smallmatrix}\right]\ }
&{\ U^{\filledstar}\oplus V}}
\text{\quad of }\underline G
\end{equation*}
(we interchanged the summands in
$V\oplus U^{\filledstar}$ of
${\mathcal M}\oplus\mathcal
 M^{\circ}$ to make it selfajoint);
 ${\mathcal M}\oplus\mathcal
 M^{\circ}$ corresponds to the
 representation
 \[
 {\mathcal
 M}^+:\ U\oplus
V^{\filledstar}
\fform \left[\begin{smallmatrix}
  0 & B^{\filledstar} \\
  A & 0 \\
\end{smallmatrix}\right]\text{\quad of }G.
\]
By \eqref{lic}, there is a set
$\ind(\underline{G})$ consisting of
the representations
 \[\mathcal
M_n(\lambda):\xymatrix@1{ {\mathbb
 C^n} \ar@<0.4ex>[r]^{J_n(\lambda
 )} \ar@<-0.4ex>[r]_{I_n} &{\mathbb
 C^n}}\qquad(\lambda \ne 0)\] and pairs of
mutually adjoint representations
 \begin{equation}\label{lhr}
 \xymatrix{
 {\mathbb C^n}
 \ar@<0.4ex>[r]^{J_n(0)}
 \ar@<-0.4ex>[r]_{I_n}
 &{\mathbb C^n}}
              \text{and}
              \xymatrix{
 {\mathbb C^n}
 \ar@<0.4ex>[r]^{I_n}
 \ar@<-0.4ex>[r]_{J_n(0)^T}
 &{\mathbb C^n}};
               \
 \xymatrix{
 {\mathbb C^{n}}
 \ar@<0.4ex>[r]^{L_n}
 \ar@<-0.4ex>[r]_{R_n}
 &{\mathbb C^{n-1}}}
               \text{and}
 \xymatrix{
 {\mathbb C^{n-1}}
 \ar@<0.4ex>[r]^{\ \ R_n^T}
 \ar@<-0.4ex>[r]_{\ \ L_n^T}
 &{\mathbb C^{n}}}
 \end{equation}
(For each matrix $M$,
$M^{\filledstar}$ is
$M^{\text{\textasteriskcentered}}$
if the fixed involution on $\mathbb
C$ is complex conjugation, or $M^T$
if the involution is the identity.)
${\cal M}_n(\lambda)$ is isomorphic
to ${\cal M}_n(\mu)^{\circ}$ if and
only if $J_n(\mu)$ is similar to
$J_n(\lambda)^{\filledstar-1}$.
${\cal M}_n(\lambda)$ is isomorphic
to a selfadjoint representation $
\xymatrix@1{ {\mathbb C^n}
 \ar@<0.4ex>[r]^{A}
 \ar@<-0.4ex>[r]_{A^{\filledstar}}
 &{\:\mathbb C^n}}$ if and only if
 $A^{\filledstar-1}A$ is similar to
 $J_n(\lambda )$. If such $A$
 exists, we fix one and denote it
 by $\sqrt[\displaystyle
 \filledstar]{J_n(\lambda)}$.

 By Fact 1, each representation of
 $G$ over $\mathbb C$ is isomorphic
 to a direct sum of representations
 of the following forms:
\begin{itemize}
  \item $\mathcal
      M_n(\lambda)^+:\ {\mathbb
      C}^n\form
      \left[\begin{smallmatrix}
  0 & I_n \\
  J_n(\lambda ) & 0 \\
\end{smallmatrix}\right]$\: if does not
exist $\sqrt[\displaystyle
 \filledstar]{J_n(\lambda)}$\ ,
 in which $J_n(\lambda)$ is
 determined up to replacement
 by $J_n(\mu)$ that is similar
 to
 $J_n(\lambda)^{\filledstar-1}$;

  \item%
 $ {\mathbb C}^n\form \:
\varepsilon\!
 \sqrt[\displaystyle
 \filledstar]{J_n(\lambda)}$\ ,
 in which $\varepsilon=\pm 1$
 if the involution on $\mathbb
 C$ is complex conjugation and
 $\varepsilon=1$ if the
 involution is the identity;

  \item $ {\mathbb C}^m \form
      J_m(0)$, which is
      isomorphic to ${\mathbb
      C}^n\form
      \left[\begin{smallmatrix}
  0 & I_n \\
  J_n(0) & 0 \\
\end{smallmatrix}\right]$ if
$m=2n$ or ${\mathbb C}^n\form
      \left[\begin{smallmatrix}
  0 & R_n^{\filledstar} \\
  L_n & 0 \\
\end{smallmatrix}\right]$ if $m=2n-1$
(these forms are obtained from
\eqref{lhr}).
\end{itemize}
The matrix $\sqrt[\displaystyle
 \filledstar]{J_n(\lambda)}$ exists
 if and only if $|\lambda |=1$ when
 the involution on $\mathbb C$ is
 complex conjugation, and $\lambda
 =(-1)^{n+1}$ when the involution
 on $\mathbb C$ is the identity.
 Respectively, one can take
 $$\sqrt[\displaystyle
 \text{\textasteriskcentered}]{J_n(\lambda)}=\sqrt{\lambda}\Delta
 _n,\qquad
 \sqrt[T]{J_n((-1)^{n+1})}=\Gamma
 _n,$$ where
\[
\Gamma_n:=\begin{bmatrix}
0&&&&\udots\\
&&&-1&\udots\\
&&1&1\\&-1&-1
\\1&1&&&0\end{bmatrix},\qquad
\Delta_n:=\begin{bmatrix}
0&&&1\\&&\udots&i
\\&1&\udots\\
1&i&&0 \end{bmatrix}
\quad\text{($n$-by-$n$)}.
\]
We obtain the following canonical
forms of a square complex matrix
$A$:
\begin{itemize}
  \item %
$A$ is *congruent to a direct
sum of matrices of the form
$$\begin{bmatrix}
  0 & I_n \\
  J_n(\lambda ) & 0 \\
\end{bmatrix},\quad
\mu\Delta_n,\quad J_n(0),$$ in
which $|\lambda|>1$ and
$|\mu|=1$;

  \item %
$A$ is congruent to a direct
sum of matrices of the form
$$\begin{bmatrix}
  0 & I_n \\
  J_n(\lambda ) & 0 \\
\end{bmatrix},\quad \Gamma_n,\quad
J_n(0),$$ in which $0\ne
\lambda\ne (-1)^{n+1}$ and
$\lambda$ is determined up to
replacement by $\lambda^{-1}$.
\end{itemize}
These direct sums are uniquely
determined by $A$, up to
permutation of summands.

\end{enumerate}
\end{example}


 \end{document}